\documentclass[10pt]{article}
\usepackage{amssymb}
\usepackage{latexsym}
\usepackage{euscript}

\newcommand{\pp}{\mathbb{P}}
\newcommand{\qq}{\mathbb{Q}}
\newcommand{\cc}{\mathbb{C}}

\newcommand{\zz}{\mathbb{Z}}
\newcommand{\sss}{\mathbb{S}}
\newcommand{\rot}{\mathbb{T}}

\DeclareFontFamily{OT1}{rsfs}{}
\DeclareFontShape{OT1}{rsfs}{n}{it}{<-> rsfs10}{}
\DeclareMathAlphabet{\curly}{OT1}{rsfs}{n}{it}

\newcommand{\cdbar}{\overline{\partial}}

\newcommand{\Pic}{\mathrm{Pic \,}}

\newcommand{\Jac}{\mathrm{Jac}}
\newcommand{\Hom}{\mathrm{Hom}}
\newcommand{\rk}{\mathrm{rk}}

\newcommand{\mgbar}{\overline{M}_g}

\newcommand{\End}{\mathrm{End}}
\newcommand{\Det}{\mathrm{Det}}

\newcounter{Universal}[section]
\renewcommand{\theUniversal}{\thesection.\arabic{Universal}}

\newenvironment{Plain}{\refstepcounter{Universal} \par \vspace{0.5cm}
\noindent {\bf (\theUniversal)}\ }{\par \vspace{0.5cm}}
\newenvironment{Italic}{\refstepcounter{Universal} \par \vspace{0.5cm}
\noindent {\bf (\theUniversal)}\ \it}{\par \vspace{0.5cm}}  
  
\newenvironment{Thm}{\begin{Italic}{\sc Theorem: }}{\end{Italic}}
\newenvironment{Prop}{\begin{Italic}{\sc Proposition:} }{\end{Italic}}

\newenvironment{Pf}{\par \noindent{\sc Proof:} }{\quad $\blacksquare$ \par
\vspace{0.5cm}}

\newenvironment{Question}{\begin{Plain} \noindent{\sc Question:
}}{\end{Plain}}

\newcounter{enum}

\newenvironment{Eqn}{\refstepcounter{Universal} $$} {\eqno \mathrm{
(\theUniversal)} $$} 

\title{Symplectic four-manifolds and conformal blocks}
\author{Ivan Smith\thanks{DPMMS, Centre for Mathematical Sciences,
Wilberforce Road, Cambridge CB3 0WB, UK.  E-mail:
i.smith@dpmms.cam.ac.uk.  MSC Classification: 53D35, 53D50, 14H60}}  
\date{}

\begin{document}
\maketitle
\thispagestyle{empty}

\begin{abstract}
\noindent We apply ideas from conformal field theory to
study symplectic four-manifolds, by using modular functors to
``linearise''  Lefschetz fibrations.  In Chern-Simons theory this
leads to the study of 
parabolic vector bundles of conformal blocks.  Motivated by the Hard
Lefschetz theorem, we show the bundles of $SU(2)$
conformal blocks associated to K{\"a}hler surfaces are Brill-Noether
special, although the 
associated flat connexions may be irreducible if the surface is simply
connected and not spin.
\end{abstract}

\section{Introduction}

This note is intended to publicise the
following juxtaposition, the potential of which is surely not fully realised.
(i) A symplectic
four-manifold $X$ can be described via Lefschetz pencils $f:X \dashrightarrow \pp^1$, which
are algebraically encoded as representations $\rho_{X,f}$ of free groups in mapping
class groups.  These representations are not canonical, but become
so (asymptotically) under a stabilisation procedure which involves a
sewing operation  
on the underlying fibres of the pencil.  (ii) Chern-Simons 
theory gives rise to (projective) linear representations $\rho_{G,k}$ of mapping class groups, once a compact Lie group $G$ and level $k$ are fixed.  These representations are not independent, but
behave coherently under sewing operations of underlying families of surfaces.
\medskip

Although the similarity above motivates our study, we are unable to take real advantage of the coherence, and accordingly achieve only modest results.  A flat connexion in a vector bundle over a punctured curve determines a parabolic bundle; the moduli space of parabolic bundles is stratification by the ``Brill-Noether'' loci of bundles which admit more holomorphic sections than predicted by Riemann-Roch.
By studying restriction maps from holomorphic bundles on a K\"ahler surface $X$ to bundles on embedded curves we will prove:

\begin{Thm}

(A) If $X$ is simply connected and not spin then $\rho_{SU(2),1}\circ\rho_{X,f}$ is irreducible.

(B) For $k \gg 0$ $\rho_{SU(2),k}\circ\rho_{X,f}$ is Brill-Noether special.
\end{Thm}

\noindent Result (A) represents the failure of a ``non-abelian'' Hard Lefschetz theorem, and is a genericity result for the $SU(2)$ Chern-Simons representations.  (The result apparently generalises from level $k=1$ to the case $k+2$ prime, but we will not prove that here; the $k=1$ case is sufficient to answer in the negative a question of Tyurin, as we discuss.)  In contrast, result (B) shows that the representations arising from K\"ahler Lefschetz pencils \emph{are} special from at least one point of view, and in principle provides
a new obstruction to integrability for symplectic four-manifolds.  (Other known obstructions that go beyond topology come from gauge theory; we contrast Result (B) with related ideas from Donaldson theory at the end of the paper.)

The next section briefly recalls background material; the third discusses the Hard Lefschetz theorem and the first result above, and the final section discusses parabolic bundles.  This paper fits into a general programme which replaces fibres of Lefschetz fibrations by moduli spaces of objects on those fibres.  Surprisingly, when the substituted moduli space is linear, the resulting object seems less tractable than e.g. when one replaces fibres by their symmetric products \cite{DS}.  Probably the right setting for these ideas has not yet been
found; we hope, despite its preliminary flavour, the paper may
encourage other people to think along these lines.  \medskip

\noindent \textbf{Acknowledgements}\footnote{The first draft of this paper was
written in Paris with the support of
  an EC Marie-Curie Fellowship HPMF-CT-2000-01013; thanks to the Ecole
  Polytechnique for hospitality and espresso, and to \emph{La Papillote}
for all the \emph{moelleux au chocolat}.}  Conversations 
with Denis Auroux, Simon Donaldson and
Graeme Segal were helpful; Andrei Tyurin's preprint
  \cite{Tyurin} has been influential \emph{passim}. 

%%%%%%%%%%%%%%%%%%%%%%%%%%%%%%%%%%%%%%%%%%%%%%%%%%%%%%%%%%%%%%%%%%%%%%%%%%%
\section{Background}

A Lefschetz pencil (or complex Morse function) on a smooth oriented four-manifold $X$ is a map $f:X\backslash\{b_1,\ldots,b_n\}\rightarrow\sss^2$ defined on the complement of a finite set, submersive away from a disjoint finite set $\{p_1,\ldots,p_{n+1}\}$, and conforming to local models $(z_1,z_2)\mapsto z_1/z_2$ near $b_j$ and $(z_1,z_2)\mapsto z_1z_2$ near $p_i$, where the $z_i$ are oriented local almost complex co-ordinates.   Donaldson \cite{Donaldson:pencils} has proved that all symplectic manifolds admit this structure. Topologically, $X$ is swept out by surfaces, finitely many of which have complex ordinary double point singularities (at the $p_i$) and all of which meet at the $b_j$ and are otherwise disjoint.  There are other helpful viewpoints: \medskip

(i) A Lefschetz pencil $f: X \dashrightarrow \sss^2$ induces a
representation $\rho_{X,f}: \mathbb{F}_n = \pi_1(\sss^2 \backslash
\{\mathrm{Crit}\}) 
\rightarrow \Gamma_g$, which is
well-defined up to global conjugation and the action of the braid
group by automorphisms of the domain.  A loop encircling one critical value maps to a positive Dehn twist about the vanishing cycle \cite{GS}. \medskip

(ii) A metric on $X$ gives a map $\sss^2\backslash\{b_j\}\rightarrow M_g$ which extends to a map $\phi_f: \sss^2\rightarrow\mgbar$ into the Deligne-Mumford moduli space of stable curves.  The homology class of the image is characterised by the fact that $\sigma(X) = 4\langle \lambda_g, [\phi_f(\sss^2)]\rangle - \delta$ where $\lambda_g$ denotes the Hodge class and $\delta$ the total number of singular fibres of the pencil.  As a consequence, one can show that $\langle \lambda, [\phi_f(\sss^2)] \rangle > 0$ always \cite{ivanhodge}.  \medskip

\noindent Donaldson's existence theorem is canonical in a certain
asymptotic limit: the closures of fibres of the Lefschetz pencil are
symplectic submanifolds Poincar{\'e} dual to $\kappa[\omega]/2\pi$, and
if the degree $\kappa$ is large enough -- depending on the particular
$X$ -- then the representation of the free group above is canonically associated to $X$
(up to global conjugation and the action of the braid group $B_n$).  There is an explicit procedure \cite{AK}, \cite{ivanmodulidivisor} which 
relates pencils of one degree $\kappa$ to pencils of a larger degree
$2\kappa$; the degree $2\kappa$ pencil is obtained
by perturbation of a degenerate family of hyperplane sections $\{s^2 + \lambda
s.t\}_{\lambda \in \pp^1}$ for $s,t$ degree $\kappa$ sections.  One feature of this stabilisation is that the
generic fibre (i.e. far away from $\lambda =0$) at degree $2\kappa$
is obtained by connect summing two
fibres of the degree $\kappa$ pencil at all their intersection points
(the base-points $\{b_j\}$ above).  This is where a family
surgery enters, in the vein of the opening remarks of the Introduction, justifying the first half of the ``juxtaposition''. \medskip

Now we briefly review Chern-Simons theory, as relevant for our needs.  Proofs and details can be found in \cite{SegalCFT},\cite{CFT}.  Fix a Riemann surface $\Sigma$ and a gauge group $G$,
for us always $SU(2), SO(3)$ or $U(1)$.  Let
$M_{2,L}(\Sigma)$ denote the moduli space of rank two stable bundles
on $\Sigma$ with fixed determinant equal to $L$; $G=SU(2)$ corresponds to $L\cong\mathcal{O}$.  The moduli space is a smooth complex
variety, closed if $deg(L)$ is odd, and with compactification
the moduli space of semistable torsion-free coherent sheaves (with 
fixed Hilbert polynomial) when $deg(L)$ is even.
The compactification locus is a copy of the Kummer
variety of the curve, which arises as the singular locus in the moduli
space if $g>2$.  The resulting projective varieties have Picard
group  $\zz$, generated by a determinant
line bundle $\mathcal{L}_{det}$ (described in more detail in the proof of
\ref{nonabeliansplits}).
Let $V_k(\Sigma)$ denote the space of holomorphic global
sections of $\mathcal{L}_{det}^{\otimes k}
\rightarrow M_{2,L}(\Sigma)$.  This is usually called the space of
\emph{conformal blocks} of level $k$ on $\Sigma$, and has dimension:  
\begin{Eqn}\label{verlindeformula}
v_k(g) \ = \ \mathrm{rank}(V_k(\Sigma_g)) \ = \ \left( \frac{k+2}{2}
\right)^{g-1} 
\sum_{j=1}^{k+1} \left( \frac{1}{\sin (j\pi/(k+2))} \right)^{2g-2}.
\end{Eqn}
(first
conjectured by Verlinde \cite{Ve}).
%; some explicit values (for
%$k=3,5$) are given in the next table.
%\begin{table}[hbt] \label{tableofranks}
%\begin{center}
%\begin{tabular}{|c|c c c c c|}
%  \hline
%$g$ & 2 & 3 & 4 & 5 & 6 \\
%\hline
%$v_3(g)$ & 20 & 120 & 800 & 5600 & 40000 \\
%$v_5(g)$ & 56 & 784 & 13328 & 241472 & Large \\
%\hline
%\end{tabular} \caption{Ranks of Verlinde bundles} \end{center} \end{table}
We suppress the group $G$, equivalently the Chern classes of the
relevant semistable sheaves, from the notation for simplicity.
The \emph{Verlinde bundle} $V_k \rightarrow M_g$ is the holomorphic vector
bundle over the moduli space of curves whose fibre at $\Sigma$ is
$V_k(\Sigma)$; that these spaces fit together to give a vector bundle follows from 
elliptic regularity. 

\begin{Thm} \cite{Hitchin,SegalCFT,Kohno} \label{flat}
Fix a gauge group $G$ and level $k$.  The holomorphic vector bundle
$V_k \rightarrow M_g$ carries a projectively flat connexion, defining
a representation $\rho_{G,k}: \Gamma_g \rightarrow \pp \End(V_k)$.
\end{Thm}

We will need several other properties of these bundles: \medskip

(i) The bundle $V_k \rightarrow M_g$ has a distinguished extension to a holomorphic vector bundle $V_k \rightarrow \mgbar$, see \cite{TUY}.  (In fact, $V_k$ has a natural parabolic structure over $\mgbar\backslash M_g$, since the flat connexion on $V_k$ provided by the theorem above has simple poles along the divisor of nodal curves, but we shall not use this.) \medskip

(ii) The first chern class $c_1(V_k \rightarrow \mgbar) = \frac{3k v_k(g)}{k+2} \lambda \in H^2(\mgbar)$.  This doesn't seem to be well-known: we give a (loop-group inspired) proof in the last section. \medskip

(iii) If $G=U(1)$ the conformal blocks are theta-functions (by lifting holomorphic sections of a line bundle over the Jacobian to periodic functions on the universal cover), cf. (\ref{abeliansplits}).  The representation $\rho_{U(1),k}$ 
factors through the automorphism group of a finite
Heisenberg group $0 \rightarrow (\zz/k\zz)^{2g} \rightarrow \widetilde{\mathrm{Aut}}
\rightarrow Sp_{2g}(\zz/k\zz) \rightarrow 0$, \cite{LB}, hence is irreducible.  \medskip

(iv) The representations $\rho_{SU(2),1}$ and $\rho_{U(1),2}$ are equal; this is part of ``rank-level duality'' \cite{CFT}, see also \cite{Kohno}, \cite{vG-dJ}. \medskip

\noindent For surfaces with marked points or boundary, there are also spaces of conformal blocks defined using moduli spaces of parabolic bundles with fixed conjugacy \cite{Kohno}.  The extended theory is a ``modular functor'' \cite{SegalCFT}, i.e. if $\Sigma = \Sigma_1\cup_{C}\Sigma_2$ where $C \subset \Sigma_i$ is a closed one-manifold, then $V_k(\Sigma) = V_k(\Sigma_1)\otimes V_k(\Sigma_2)$.
Modularity implies that 
the Verlinde bundles are compatible under sewing in the following
sense: given two families of Riemann surfaces over a base $B$ and
identification diffeomorphisms of (some subsets of) the boundaries
over the base, the projectively flat connexion in the vector bundle
associated to the glued surfaces splits in the tensor product
decomposition of the vector bundle. This justifies the second half of the ``juxtaposition'' from the Introduction.

%%%%%%%%%%%%%%%%%%%%%%%%%%%%%%%%%%%%%%%%%%%%%%%%%%%%%%%%%%%%%%%%%%%%%
\section{Hard Lefschetz}

Let $f:X \rightarrow \sss^2$ be a
Lefschetz fibration; then $H^1(X)$ is the subgroup of monodromy
invariants of $H^1(F)$, and if $X$ is K\"ahler the Hard Lefschetz theorem asserts that $H^1(X) \subset
H^1(F)$ is a \emph{symplectic} subspace, equipped with the
non-degenerate skew-form $(a,b) \ = \ \int_X \omega \wedge a \wedge b$.
(For general symplectic pencils the monodromy representation
will not be completely reducible and the invariant subgroup of
$H^1(F)$ and the subgroup generated by Poincar{\'e} duals of vanishing
cycles will 
not be orthogonal with respect to cup-product.) 

\begin{Prop}\label{abeliansplits}
If $X$ is K{\"a}hler and $b_1(X)>0$ then the
representations $\rho_{U(1),k} \circ \rho_{X,f}$ are reducible for every $k$.
\end{Prop}

\begin{Pf}
The abelian Verlinde representations
$\rho_{U(1),k}$ factor through the metaplectic representation of the
double cover $Mp_{2g}(\zz)$ of the symplectic group.  The vector space
for a genus $g$ surface at level $k$ is a $k^g$ dimensional space of
theta-functions:

$$H^0(\Jac(\Sigma_g)~; \mathcal{L}^k) \ = \ \langle \vartheta_i^j \ |
\ 1\leq i \leq k, \, 1\leq j\leq g \rangle$$

$$ \vartheta_i^j (\underline{z},\tau) \ = \ \sum_{\underline{n} \in
  \zz^g + e^j_i} \exp \left( \frac{i\pi}{k} (\underline{n}^t \cdot
  \tau \cdot \underline{n}) + 2i\pi (\underline{n}\cdot \underline{z})
\right)$$ 

\noindent where $\tau \in \frak{h}_g$, the Siegel upper half-space,
$\underline{z} \in \cc^g$ and $e_i^j$ has value $i/k$ in the $j$-th
position and zeroes elsewhere.  This set of generators for the space
of conformal blocks leads to a ``factorisation'' property.

Suppose a symplectic vector space $W$ is written as a product of
symplectic subspaces $W = U \oplus U'$ of dimensions $2h, 2g-2h$
respectively.  We have a natural inclusion of symplectic groups
$Sp_{2h}(\zz) \times Sp_{2h-2h}(\zz) \hookrightarrow Sp_{2g}(\zz)$.
On the image of this inclusion, the metaplectic representations
factorise as a tensor product:

$$\rho_{U(1),k}[g](A\oplus B) \ = \ \rho_{U(1),k}[h](A) \otimes
\rho_{U(1),k}[g-h](B)$$ 

\noindent in an obvious notation.  This is proved by the following
calculation with exponentials.  Fix a vector $\underline{l} =
\underline{l}_1 \oplus \underline{l}_2$ which indexes a particular
choice of $\vartheta$-characteristic \cite{LB} which above is given by
the choice of labels $e_i^j$.  Then we have

$$\sum_{\underline{m} \in \zz^{h+(g-h)}} \exp \left \{ \frac{i\pi}{k}
  \left( (\underline{m}+\underline{l})^t \left( \begin{array}{cc} A &
        0 \\ 0 & B \end{array} \right) (\underline{m} + \underline{l})
  \right) \right \}$$

$$= \ \sum_{\underline{\alpha} \in \zz^h, \, \underline{\beta} \in
  \zz^{g-h}} \exp \left \{ \frac{i\pi}{k} \big(
  (\underline{\alpha}+\underline{l}_1)^t A (\underline{\alpha} +
  \underline{l}_1) + (\underline{\beta} + \underline{l}_2)^t B
  (\underline{\beta} + \underline{l}_2) \big) \right \} \ = $$

$$ \left( \sum_{\underline{\alpha} \in \zz^h} \exp
  \left\{ \frac{i\pi}{k} \big( (\underline{\alpha}+\underline{l}_1)^t
    A (\underline{\alpha} + \underline{l}_1) \big) \right\} \right)
\left( \sum_{\underline{\beta} \in \zz^{g-h}} \exp
  \left\{ \frac{i\pi}{k} \big( (\underline{\beta}+\underline{l}_2)^t
    B (\underline{\beta} + \underline{l}_2) \big) \right\} \right).$$

\noindent For a K{\"a}hler
Lefschetz fibration there is a symplectic splitting $H^1(F) =
H^1(X) \oplus \mathrm{Ann}(H^1(X))$, and the homological monodromy is
trivial on the first factor.  It follows that the Verlinde
representation is of the form $(\mathrm{id} \otimes \phi)$ which is
clearly reducible.
\end{Pf}

\noindent To generalise, we think of the reducibility arising from the presence of line bundles on the surface ($b_1(X)>0$); Gieseker and O'Grady \cite{oG} have shown that every projective surface has a positive-dimensional moduli space of stable bundles with trivial determinant and fixed $c_2=r$ once $r$ is sufficiently large.
However, there is no
analogous reducibility property for the representations in general.
To see this we begin with a
statement about the monodromy groups of Lefschetz
pencils, obtained jointly with Denis Auroux.

\begin{Prop}
Let $X$ be a symplectic manifold.  If $X$ is spin, the homological
monodromy representation of any Lefschetz pencil is not onto
$Sp_{2g}(\zz)$.  If $X$ is not spin and $H_1(X;\zz/2\zz)=0$, the homological monodromy
representation is surjective for any pencil given by stabilisation
(degree doubling).
\end{Prop}

\begin{Pf}
The key ingredient is a result of Janssen \cite{Ja} which in turn
relies on work of Gabrielov and Chmutov \cite{Ch}; we give their
results translated 
into our language. Take a pencil of curves satisfying the constraints
that (i) two vanishing cycles have homological intersection number one
and (ii) all the (homology classes of the) vanishing cycles are
conjugate under the monodromy 
group of the pencil.  There are pencils that violate these
conditions, for instance the genus two pencil on $\rot^2 \times
\sss^2$.  However, the first condition always holds after
stabilisation, by a quick look at the pictures of Auroux and Katzarkov
\cite{AK}, whilst the second condition holds for large enough degree
by a result of Amoros, Munoz and Presas \cite{AMP}.  (In the algebraic setting,
the second condition follows from the
irreducibility of the dual variety for a projective embedding of the
surface, together with the Lefschetz hyperplane theorem; these imply
that the fundamental group of the complement of the dual variety is
normally generated by one element.)   Chmutov proves that for any such
pencil of curves, the homological monodromy contains the kernel of the
natural map $Sp_{2g}(\zz) \rightarrow Sp_{2g}(\zz/2\zz)$ (this is
generated by the squares of Dehn twists).  Janssen builds on this to
deduce that either the monodromy of a pencil is contained in the
hyperelliptic mapping class group (which is an easy exceptional case,
not preserved by doubling), or is full, or maps onto the subgroup of
$Sp_{2g}(\zz/2\zz)$ which preserves a quadratic form. 

Now we remark that a four-manifold with $H_1(X;\zz/2\zz)=0$ is spin iff the
intersection form is even.  Using this, it is easy to see that $X$ is
spin iff the associated Lefschetz fibration over a disc (blow up the
base points and remove a smooth fibre) is spin.  However, as Stipsicz
points out in \cite{Stipsicz3}, this Lefschetz fibration has a
distinguished handle decomposition, in which handles are added to
$\Sigma \times D$ along the vanishing cycles with framing $-1$.  A
spin structure $q: H_1(\Sigma, \zz/2\zz) \rightarrow \zz/2\zz =
\{0,1\}$ on $\Sigma \times D$ 
extends across a handle whenever the associated vanishing cycle
evaluates to $+1$.

There are two quadratic forms $q$ on $H_1(\Sigma_g, \zz/2\zz)$, to
isomorphism, determined by their Arf
invariant (in suitable symplectic bases they correspond to $(x,y)
\mapsto \sum x_iy_i$ or to $(x,y) \mapsto \sum x_iy_i + x_1^2 + x_2^2$).
If either of these forms is preserved by the monodromy of a pencil, $X$
is spin; fixing a form
fixes either an even or odd spin structure on a fibre of the pencil,
hence a spin structure on $D \times \Sigma_g$, and this extends over
the handles added along the vanishing cycles.  (Conversely, a spin
structure on $X$ induces one on the codimension zero subset $D \times
(\Sigma_g \backslash \{b_j\})$ and hence one on $\Sigma$, the parity of
which is reflected in the monodromy group of the pencil.) 
\end{Pf}

\noindent It is possible that some symplectic manifolds admit
Lefschetz pencils of arbitrarily high degree whose monodromy group is
the entire mapping class group; the amazing computations of Auroux and
Katzarkov 
\cite{AK} give an obvious route to attack such a question. Roberts \cite{Ro} has shown that the
$SU(2)$-theory mapping class group representations on
conformal blocks remain irreducible whenever $k+2$ is prime.
His combinatorial method of proof, together with the pictures of
\cite{AK}, strongly suggests the  
composite representations $\rho_{SU(2),k}\circ \rho_{X,f}$ are
irreducible for all such $k$. 
The previous Proposition, together with Remarks (iii) and (iv) after (\ref{flat}), prove Theorem (A) from the Introduction.  \medskip

The Proposition is also relevant to a question of Tyurin; to see this we need to gather some facts about the behaviour of rank two
stable bundles under restriction.  We will use
a strong form of the restriction theorem that can be garnered from
\cite{Tyurin}, \cite{Fr} and \cite{shaves}.
Let $X$ be a projective surface with hyperplane class $H$.  Let $C$
denote an arbitrary smooth element of the linear system $|NH|$.
\begin{itemize}
\item Let $V$ be an $H$-stable $SU(2)$ or $SO(3)$ bundle on $X$. If $N
  > -p_1(ad V)$ then $V$ restricts to a stable bundle on $C$. 
\item For $N \geq N(p_1) \gg 0$ restriction defines a holomorphic
  embedding $M_r(X) \rightarrow M_{2,\mathcal{O}}(C)$ for each moduli
  space with $0\leq r \leq -p_1$.  
\end{itemize}
For non-hyperelliptic curves the determinant line is very ample, i.e. $M_{2,\mathcal{O}}(C) \hookrightarrow \pp(V_k(C)^*)$ for all $k\geq 1$, and we can compose the restriction maps with these embeddings.
\begin{Question}(Tyurin,\cite{Tyurin},\cite{Tyurin3}) \label{tyurinconj}
With notation as above, at level $k=1$, does the canonical connexion
preserve the image of the  
fibrewise restriction map as we vary $C$ amongst smooth curves (all or none of which are hyperelliptic) in its linear system? 
\end{Question}

\noindent Tyurin conjectured this was true as we vary $C$ locally and
asked when it was true globally; in which case we'll say that the stable bundles
on the surface $X$ are \emph{globally invariant}.  Not surprisingly,
the irreducibility property for $\rho_{SU(2),1}\circ\rho_{X,f}$
described above implies that global invariance is exceptional.  We'll call a polarisation $H$ of an algebraic surface $X$ ``appropriate'' if (i) $H$ is even and $K_X$ is even or (ii) $H=2n\tilde{H}+K_X$ for some arbitrary polarisation $\tilde{H}$ and $n>0$, when $K_X$ is odd.  These obviously exist, so we can consider ``appropriate'' pencils $\pp^1 \subset |\kappa H|$.

\begin{Prop}\label{nonabeliansplits}
Let $X$ be a non-spin surface of general type with
$H_1(X;\zz/2\zz)=0$; fix an appropriate polarisation $H$ on $X$ and an appropriate pencil of degree $\kappa$. The stable bundles on $X$ of Chern class $c_2$ are not globally invariant for $\kappa$ sufficiently large.
\end{Prop}

\begin{Pf}
Complex surfaces have non-empty moduli spaces $M_r(X)$ of $SU(2)$ bundles (indeed the moduli space contains smooth points as soon as $c_2=r > b_+(X)+2$ \cite{Fr}), with well-understood compactifications.  We will
take the second Chern class large enough for the moduli space to be``generically smooth'' \cite{Donaldson:polynomial}, meaning it is of the expected dimension and the singular locus has codimension greater than one (explicit bounds on the required $c_2$ are known).  The singularities are normal and algebraic sections of line
bundles over the 
moduli space are uniquely determined by their behaviour
on the open smooth locus \cite{shaves}.  With this background, the
following calculation will be formal but can be
put on a solid footing in a familiar way.

Fix a smooth curve $C \in |\kappa H|$ in $X$, $\kappa \gg 2r$; by restriction we can suppose $M_r(X) \subset M_{2,\mathcal{O}}(C)$.   Let $L_{\kappa} =
\mathcal{L}_{det}|_{im(r)} \rightarrow M_{stab}(X)$ be the line bundle
on $M_r(X)$ given by restricting the determinant line from $M_{2,\mathcal{O}}(C)$.  We claim that $h^0(L_{\kappa})$ grows at most polynomially with
$\kappa$.  This is not completely trivial, since $\Pic(M_r(X))$
is complicated and the dependence of $L_{\kappa}$ on $\kappa$ is not
linear \cite{shaves}, but it reduces to an argument of Donaldson from
\cite{Donaldson:polynomial}.
The determinant line bundle on
$M_{stab}(C)$ at a bundle $E \rightarrow C$ is defined as having fibre 
\begin{Eqn} \label{determinant}
(\mathcal{L}_{det})_E \ = \ \Lambda ^{max}H^0(E \otimes K_C^{1/2})
\otimes \Lambda ^{max} H^1 (E^* \otimes K_C^{-1/2})^{-1};
\end{Eqn}
by Serre duality (and recalling that $E^* \cong E$ for $SU(2)$ bundles), we can simplify this to $(\mathcal{L}_{det})_E \ = \
\Lambda ^{max}H^0(E \otimes K_C^{1/2})^2$.  The choice of spin
structure on $C$ plays no serious role, as explained in \cite{DK}, p.382.
The divisor sequence on $X$ 
$$0 \rightarrow F \rightarrow F(D) \rightarrow F(D)|_D \rightarrow 0$$
with $F=E\otimes \mathcal{O}(-NH+K_X)^{1/2}$ and $D=C\in |NH|$,
and the adjunction formula $K_C = (K_X +
\mathcal{O}_X(C))|_C$, gives the long exact sequence in cohomology
$$0 \rightarrow H^0(E\otimes \mathcal{O}(-NH+K_X)^{1/2}) \rightarrow
H^0(E\otimes (NH+K_X)^{1/2}) \rightarrow \cdots $$
$$ \cdots H^0 (E|_C \otimes K_C^{1/2})
\rightarrow H^1(E\otimes \mathcal{O}(-NH+K_X)^{1/2}) \rightarrow \cdots$$
For $H$ ample and $N$ large enough the
bundles $F\otimes \mathcal{O}(-NH+K_X)^{1/2}$ will have no cohomology
except in the top dimension.  Taking
determinants, we have that 
$$(\mathcal{L}_{det})_{E|_C} \ = \ \Lambda ^{max}H^0(E \otimes
(NH+K_X)^{1/2})^2$$ 
where the relevant bundles and cohomology group all live on $X$.
Although the dependence of this bundle on $N$ is not linear as we vary
the choice of linear system, over the open smooth locus of
$M_r(X)$ where a universal bundle on $X \times M_{stab}(X)$
exists, all of the bundles $L_{\kappa}$ are pushforwards
$det(R^0 \pi_* \mathbb{E}\otimes L^{\kappa})$, where $\mathbb{E}$ is
the universal bundle, $L$ is the bundle on $X$ with first Chern class
$H$ pulled back to $X \times M_r(X)$, and $\pi$ denotes the
projection to $M_r(X)$.  Moreover, once $N$ is large, Li has shown the bundles
$L_{\kappa}$ are ample \cite{JunLi},
\cite{shaves}.  It follows that $h^0 (L_{\kappa})$ can be computed as an
Euler characteristic by Grothendieck-Riemann-Roch, hence is a polynomial in the characteristic
classes of $\mathbb{E} \otimes L^{\kappa}$ and the Pontrjagin classes
of $M_r(X)$, yielding the claim. \medskip

\noindent With this established, suppose
stable bundles on $X$ are globally invariant for a Lefschetz pencil $f$ of
high degree.  There is a monodromy representation $\rho_{SU(2),1}
\circ \rho_{X,f}~: \mathbb{F}_n \rightarrow \pp GL(N_{\kappa})$.  Each matrix
in the image of this representation preserves the subvariety
$M_{stab}(X) \subset \pp^{N_{\kappa}}$ and hence in particular preserves the
locus of hyperplanes 
$$\{ h \in (\pp^{N_{\kappa}})^* \, | \, M_{stab}(X)
\subset h \},$$
that is the set of conformal blocks vanishing completely on the subset
of restriction.  These hyperplanes are exactly the rays in the kernel
of the natural map from $H^0(V_1(C)) \rightarrow H^0(L_{\kappa})$.  The rank
of the first group grows exponentially with the degree $\kappa$ of the
pencil on $X$, whereas the rank of the second group grows polynomially
by the above.  Hence, once $\kappa$ is large enough, the kernel is
non-empty.  It is obviously not full; the representation
therefore admits a non-trivial  invariant subspace, contradicting
Theorem (A). 
\end{Pf}

\noindent Global restrictions on monodromy shed no light on the ``local'' version of Tyurin's question, which in any case can apparently not be sensibly formulated for symplectic as opposed to K\"ahler Lefschetz pencils.  However, there are properties which make sense in
general and which appear special in the K\"ahler context, which we
address next.

%%%%%%%%%%%%%%%%%%%%%%%%%%%%%%%%%%%%%%%%%%%%%%%%%%%%%%%%%%%%%%%%%%%%%%%5
\section{Brill-Noether}

For every fibre genus $g$, level $k$ and number of critical
fibres $r$, fix once and for all a model of the symplectic
representation space 
$$\Hom_+(g,k,r;c_1) \ = \ \Hom_+(\pi_1(\sss^2\backslash
\{p_1,\ldots,p_r\}), \pp U_{v_k(g)} \big/ \langle \mathrm{Conj}\rangle.$$
Here the subscript $+$ indicates that we are fixing the holonomy data
at each puncture to conform to the matrix $\rho_k(\tau)$ given by a
positive Dehn 
twist in a non-separating curve, and
$c_1$ denotes the topological degree of the bundles.  This space is a
connected symplectic orbifold, whose dimension is given by 
$$\dim (\Hom_+(g,k,r)) \ = \ r\dim(C_{\tau}) - 2\dim(\pp U_{v_k(g)})$$
where $C_{\tau}$ denotes the conjugacy class in $\pp U_{v_k(g)}$ of
the matrix $\rho_k(\tau)$.  If we fix a projective unitary representation of the mapping class
group, a Lefschetz fibration gives rise to a point of the quotient of $\Hom_+$
by the braid group action of Hurwitz moves.  The braid group acts ergodically
\cite{Goldman},
so in the absence of invariant open sets one can look for
invariant stratifications.
\medskip

\noindent Any
choice of complex structure $J \in M_{0,r}$ on $\sss^2 \backslash
\{p_i\}$ defines 
a projective moduli space of parabolic bundles $\mathcal{M}_{par}(J)$,
with the flags and
monodromy at each puncture fixed to be the same local model, and a
homeomorphism $\psi_J: \mathcal{M}_{par}(J) \stackrel{\sim}{\longrightarrow}
\Hom_+$.  The space of parabolic bundles $E$ carries a natural
stratification, given by the upper semicontinuous function $E \mapsto
h^0(E^*)$ in the case where $c_1 > \mathrm{rank}(E)$; the union of all
the lower strata $\{ E \, | \, h^0(E^*)>0 \}$ is the
\emph{Brill-Noether} locus, a complex subscheme which is
carried by $\psi_J$ to a closed real subvariety of $\Hom_+(g,k,r;c_1)$.
Taking the union over the images of these subvarieties as we vary $J
\in M_{0,r}$ defines a sequence of braid-group invariant
subsets of $\Hom_+$.  Each of these 
is nowhere dense for large $k$:
the Brill-Noether loci may have excess dimension, but
their actual codimension grows with the
virtual codimension, hence with $k$.  By contrast, the 
space $M_{0,r}$ is
a smooth complex manifold of dimension $r-3$ independent of $k$.  
The dense open subset 
$$\mathcal{U} \ = \ \{ \rho \in \Hom_+ \ | \ h^0(\rho^*;j) = 0 \
\mathrm{for \ every} \ j \in M_{0,n} \}.$$
comprises the parabolic bundles
which are Brill-Noether general for \emph{every} complex structure on the base (there are moduli here since there are marked points).  Note that the condition $c_1(E) > \mathrm{rk}(E)$ holds for conformal block bundles over $\pp^1$ by Remark (ii) after (\ref{flat}), proved below. \medskip

\noindent  The Arakelov-Parsin theorem shows that in fact only finitely many conjugacy classes of representation $\rho$ are realised by K\"ahler pencils, but it seems very hard to characterise or identify properties of this distinguished finite set of braid-group orbits, which lends the following some interest.  Fix $c_2=r$ large enough for moduli spaces of bundles on $X$ to be singular only in high codimension and fix a Lefschetz pencil of degree $\kappa \gg r$ so that restriction maps are well-defined embeddings.

\begin{Thm}\label{brillnoether}
If $X$ is K{\"a}hler and $f$ is a pencil of degree $\kappa$, then $\rho_{SU(2),k}
\circ \rho_{X,f} \in \Hom_+\backslash \mathcal{U}$ for all sufficiently large levels $k$.
\end{Thm}

\begin{Pf} If $C \in |\kappa H|$, the
\emph{restriction kernel} is the subspace of
conformal blocks $\{ s\in V_k(C) \ | \ s|_{M_r(X)} \equiv 0 \}$ which
vanish identically on the image of the restriction map on $M_r(X)$.
If $k$ is sufficiently large, as we vary $C$ in its linear system
we claim there is a short exact sequence 
$$0 \rightarrow RK_k \rightarrow V \rightarrow \underline{\cc}^a
\rightarrow 0$$
of vector bundles over the base $\pp^1$.
The restriction of the
determinant line bundle $\mathcal{L}_{det}$ 
from $M_{2,\mathcal{O}}(C)$ to $M_r(X)$ gives an ample line bundle
$\mathcal{L}_{rest}$ on the latter
space which Tyurin proves does not depend on the choice of (smooth or nodal) curve
$C$ within its linear system \cite{Tyurin}, and we need to see that all sections of this line bundle extend to $M_{2,\mathcal{O}}(C)$.
The cokernel of the restriction map is given by
$H^1 (\mathcal{I}_{M_r(X)} \otimes \mathcal{L}_{rest}^{k})$, where
$\mathcal{I}$ is the ideal sheaf.   Since $\mathcal{L}_{rest}$ is ample,
this higher cohomology group will eventually vanish. This gives the claim. \medskip

\noindent For a Lefschetz fibration $f~: X \rightarrow \sss^2$ the index of
the d-bar operator on the dual of the Verlinde bundle is negative (cf. the start of this section).
Using Grothendieck's theorem that all vector bundles on the line split, together with (\ref{firstchern}), the most stable splitting type for the bundle $V_k^*$ will be $\mathcal{O}(1-3\lambda) \oplus \cdots \oplus
\mathcal{O}(1-3\lambda) \oplus \mathcal{O}(-3\lambda) \oplus \cdots \oplus \mathcal{O}(-3\lambda)$, where the ratio of the number of factors of the first sort to the total rank goes to zero as $k$ increases.  (In other words, the index is more negative than the rank, and the ratio of $c_1$ to rank as $k\rightarrow \infty$ approaches $-3\lambda$ from above.)  Since $\langle \lambda, [\sss^2]\rangle  >0$ we see that generically (in the space of holomorphic structures) $V_k^*$ has no sections; but restriction kernels give rise to sections of $V_k^*$.
\end{Pf}

We still owe the computation of $c_1(V_k)$.  This requires a background remark on the
conformal field theory connexions.  Fix a level $k$, and consider the
family of connexions in bundles over moduli spaces $M_g^n$ (of genus
$g$ curves with $n>0$ parametrised boundary components) that arise in
$SU(2)$ Chern-Simons theory. Each connexion has scalar curvature $c$, and the scalar -- normalised with respect to a natural K\"ahler form -- is \emph{independent} of $g$ and $n$ (determined only by the level $k$, in fact $c=3k/(k+2)$).  This is proved in Segal's loop group framework \cite{SegalCFT}, by using the modularity of the Verlinde bundles to ``localise'' the curvature computation to the case where the surface is an annulus.
There is a
determinant line bundle over the moduli space $\mgbar$, with fibre
given by the determinant line of the underlying $\cdbar$-operator on
the surface: $L_{\mathrm{det}}(C) \ = \ \Lambda^g H^0(K_C)^*$, and 
with  $c_1(L_{\mathrm{det}}) = -\lambda$.  Segal observed in
(\cite{SegalCFT}, Appendix B) that determinant lines give modular
functors (in particular $L_{det}(\Sigma) = L_{det}(\Sigma_0)$
canonically, when $\Sigma_0 = \Sigma \backslash D$), and that the scalar curvature of the associated theory is $c=-2$.

Recall that $\mgbar$ is a homology manifold,
and by \cite{Wolpert} there is an integral basis for $H^2(\mgbar)$
comprising the 
classes $\lambda$ and the Poincar{\'e} duals of the
subvarieties defined by curves which are separated by a node into
components of genus $i$ and $g-i$, with $0\leq i \leq [g/2]$. 

\begin{Prop} \label{firstchern}
The first Chern class $c_1(V_k) = \frac{3k}{k+2} \rk(V_k) \lambda \in
H^2(\mgbar)$. 
\end{Prop}

\begin{Pf}
The tensor
product of modular functors is also a modular functor with the central
charge behaving additively, and it follows that if we take the
association 
\begin{Eqn} \label{tensorflat}
\Sigma \ \mapsto \ \ V_k(\Sigma)^{(k+2)} \otimes \Det(\Sigma)^{3k}
\end{Eqn}
we define a \emph{flat} (not just projectively flat) vector
bundle over the moduli space $M_g$.  Hence this bundle has trivial
first Chern class, from which we quickly deduce that the given formula
holds in $H^2(M_g)$.  Strictly speaking, this argument only applies as
it stands to spaces $M_g^n$ with $n>0$, because the ``loop group definition'' of the connexions is only valid for surfaces with boundary (for a closed surface $\Sigma$ Segal defines $V_k(\Sigma)=V_k(\Sigma\backslash D)$ and shows this is independent of the disc $D$ which is removed).  However, as will become clear
below, $H_2(M_g)$ is generated by a 2-cycle which admits a lift to
$M_g^1$ (arising from a fibred four-manifold $X$ which has a section $D$).
Removing a neighbourhood of $D\subset X$ gives a \emph{coherent}
family of decompositions of the fibres into open surfaces union discs,
and enables one to reduce the computations for closed surfaces to
surfaces with non-empty boundary.  This shows the required formula does hold in $H^2(M_g)$; to lift this to $H^2(\mgbar)$ we proceed as follows. \medskip

\noindent Fix $g>2$ and construct a surface bundle with fibre
genus $g-1$, with two disjoint sections of non-zero square $-l$, and with
total space having non-zero signature $a$.  (This can be easily done by modifying appropriate Lefschetz fibrations, as in \cite{KKS}.)  Gluing the sections
together and picking a fibrewise metric gives a family of nodal curves
of genus $g$, and the base defines a curve $C \subset
\mgbar$, the image of a curve $\tilde{C} \subset M_{g-1}^2$.  $C$ lies entirely in the divisor $\Delta$
of nodal stable curves; the normal bundle $\nu_{\Delta/\mgbar}$ is
canonically identified along $C$ with the tensor product of the
tangent spaces to the two exceptional sections.  It follows that
$C\cdot \Delta = -2l$.  Also fix a smooth surface bundle of genus $g$,  
parametrised by a curve $\Sigma_{sm} \subset M_g$, with total space
having non-zero signature: $\Sigma_{sm} \cdot \lambda_g = N \neq
0$. We can assume $\Sigma_{sm}$ is the image of a surface $\tilde{\Sigma}_{sm} \subset M_g^1$ by insisting the surface bundle has a section. \medskip

\noindent  Given $\Sigma \subset \mgbar$, a curve arising
from a fibration of curves with no reducible members, the homology
class of $\Sigma$ is completely determined by the numbers $\Sigma
\cdot \Delta = t, \Sigma \cdot \lambda = t'$, since
$H^2(\mgbar, \zz) = \langle \lambda, PD[\Delta_i]
\rangle$, with $\Delta_i$ the components of the divisor of stable
curves.  Trivial algebra yields that
$$[\Sigma] \ = \ \frac{1}{N} (t+\frac{t'a}{2l})[\Sigma_{sm}] -
\frac{t'}{2l}[C]$$ 
in $H_2(\mgbar; \qq)$.  
Linearity and modularity now reduce the computation of $\langle c_1(V_k),[\Sigma]\rangle$ to the analagous pairings with $\tilde{\Sigma}_{sm}$ and $\tilde{C}$, which can be computed in terms of pairings with the relevant determinant lines by our initial reduction.  From \cite{moduliofcurves}, \cite{SegalCFT}, the Hodge class $\lambda_g$ pulls back under the
obvious map $\overline{M}_{g-1}^2\rightarrow\overline{M}_{g-1,2}\rightarrow\Delta$ to the Hodge
class $\lambda_{g-1}$ on $\overline{M}_{g-1}$ (lifted via the
forgetful map to the 
moduli space of curves with two marked points) so $C \cdot
\lambda_g = a$. 
The result now follows from (\ref{tensorflat}). 
\end{Pf}

\noindent Since $V_k \rightarrow M_g$ carries a projectively flat connexion, the higher Chern classes are incidentally given by $c_i = {n\choose i} (\frac{c_1}{n})^i$, where $n=v_k(g)$ is the rank.  This completes our treatment of Theorem (B); one can go a little further, and compare it with the kind of information obtained using Donaldson invariants. We will just outline the connection.
Let $(X,\omega)$ be an integral symplectic manifold, and fix a level
$k$ and a Lefschetz pencil $f$.  The
$k$-depth of $(X,f)$ is
defined as 
$\max_{j\in M_{0,r}} \{ h^0(V_k(f_N)^*;j) \}$;
since we have factored out the choice of $j$, this is a symplectic
invariant of the Lefschetz pencil.  It refines Theorem (B) in the sense that it measures the depth of the braid group orbit of $\rho_{SU(2),k}\circ\rho_{X,f}$ in the Brill-Noether stratification.  \medskip

There is a homeomorphism \cite{Donaldson:polynomial}
between (i) the space of instantons on an $SU(2)$ bundle $E
\rightarrow X$ with 
  $c_2(E)=r$, where $X$ is equipped with its K{\"a}hler metric, and (ii)
the space of stable holomorphic bundles topologically equivalent
  to $E$.  Donaldson invariants are suitable intersection pairings on
  $M_r(X)$, while the depth invariants are related to
  $h^0(\mathcal{L}_{rest}^k)$, which for large $k$ is also an intersection
pairing by Riemann-Roch.  This pairing gives a bound on the size of trivial quotients of $V_k$ for one 
holomorphic structure on the base of the Lefschetz fibration, hence an
estimate on the supremum over all complex structures, and hence a lower bound on the depth. 

%%%%%%%%%%%%%%%%%%%%%%%%%%%%%%%%%%%%%%%%%%%%%%%%%%%%%%%%%%%%%%%%%%%%%%%%%%%
%%%%%%%%%%%%%%%%%%%%%%%%%%%%%%%%%%%%%%%%%%%%%%%%%%%%%%%%%%%%%%%%%%%%%%%
\begin{small}

\providecommand{\bysame}{\leavevmode\hbox to3em{\hrulefill}\thinspace}

%
%\bibliographystyle{amsplain}
%\bibliography{main}
\end{small}

\end{document}